\documentclass{article}

\usepackage{amssymb}
\usepackage{dsfont}

\newtheorem{lem}{Lemma}
\newtheorem{defi}[lem]{Definition}
\newtheorem{theo}[lem]{Theorem}

\newtheorem{rem}[lem]{Remark}

\newtheorem{fait}[lem]{Fact}

\newcommand{\proof}{\noindent{\bf Proof.}~}
\newcommand{\qed}{\ \hfill$\square$\bigskip}

\title{STABILITY AND BOUNDED BALLS OF FREE PRODUCTS}
\author{Azadeh Neman}
\date{\today}

\begin{document}
\maketitle

\section{Motivation and result}

In a series of papers starting in \cite{MR1863735} and 
culminating in \cite{Sela2007}, Z. Sela proved that 
free groups, and more generally torsion-free hyperbolic groups, 
have a stable first-order theory. The question of the stability of the free product of 
two arbitrary stable groups has then been raised by E. Jaligot with, seemingly, 
the reasonable conjecture of a positive answer \cite{JaligotPekin}. 
However, a full answer seems to become a very large project of generalization, 
from free groups to free products, of the above mentioned work. 
The first step in this process is the description of varieties, that is the understanding of 
Makanin-Razborov diagrams in free products. This tranfer from \cite{MR1863735} 
(and \cite[\S1]{SelaHyp}) is represented in the work in progress \cite{JaligotSela08}. 

In the meantime, we provide here a very preliminary --- or somehow experimental --- result in the direction of the stability of free products of stable groups, 
restricting ourselves to quantifer-free definable sets and to bounded balls of free products. Notice that we consider here a {\em fixed} group, and not a class of groups as in \cite{MR2395048} or \cite{JaligotSela08}. 

Let $G\ast H$ be the free product of two groups $G$ and $H$. 
Any element of $G\ast H$ has a unique representation in normal form, i.e. it is a word in letters, that is elements of $G$ and $H$ \cite{LyndonSchupp77}. We adopt the convention that the identity element of $G\ast H$ is represented by the identity of $G$ (and not of $H$), so that the representation is unique. The notion of length is then defined in the obvious sense for these uniquely expressed normal forms of elements of $G\ast H$. 

For an integer $r\geq 1$, let $B_{r}(G\ast H)$ denote the ball of radius $r$, which is the set of elements of $G\ast H$ of length $\leq r$. Our experimental result is the following. 

\begin{theo}\label{MainTheo}
Let $G\ast H$ be the free product of two qf-stable groups $G$ and $H$. 
Let $w(\overline{x},\overline{y})$ be a group word and let $r\geq 1$ be a natural number. Then there exists a natural number $n$ (depending on $w$ and, a priori, on $r$) 
bounding the set of natural numbers $m$ for which there exists 
$\overline{a}_{1}$, ..., $\overline{a}_{i}$, ..., $\overline{a}_{m}$ and 
$\overline{b}_{1}$, ..., $\overline{b}_{j}$, ..., $\overline{b}_{m}$ 
in $B_{r}(G\ast H)$ such that $w(\overline{a_{i}},\overline{b_{j}})=1$ if and only if $i\leq j$. 
\end{theo}

Nonstandard notions will be defined shortly. Since we work in a fixed model, 
we cannot use the Infinite Ramsey Theorem just by compactness, as usual in model theory. 
Hence the present proof of Theorem \ref{MainTheo} mostly uses the 
Finite Ramsey Theorem (see \cite[Remark 0.2.7]{Wagner(Book)97}).  

\begin{fait}[Finite Ramsey Theorem]\label{finRam}
For every triple $(k,n,m)$ of natural numbers there is some natural number 
$R(k,n,m)$ such that whenever unordered $n$-tuples of a set of size at least 
$R(k,n,m)$ are painted in $k$ colors, then there is a monochromatic subset of size $m$. 
\end{fait}

\section{Technicalities and proofs}

For the basic notions of first-order logic and model theory we refer to \cite{Hodges(book)93}. As in \cite{MR2395048}, we need versions of the combinatorial properties of \cite{Shelah90} not related to a complete theory. 

\begin{defi}
Let $\cal M$ be a structure and $\phi(\overline{x},\overline{y})$ a 
formula in the same language (possibly with parameters in the domain $M$ of $\cal M$). 
Let $B$ be a subset of $M$. We say that $\phi$ is {\em stable relative to} 
$(\mathcal{M}, B)$ if there exists a maximal finite $m$ for which there 
exists $\overline{a}_{i}$ and $\overline{b}_{j}$ in $B$, $1\leq i,j \leq m$, 
such that $\phi(\overline{a}_{i},\overline{b}_{j})$ holds in $\cal M$ 
if and only if $i\leq j$. 
The maximal such $m$ is then called the {\em stability index} of $\phi$ 
relative to $(\mathcal{M}, B)$ 
\end{defi}

\begin{fait}\label{FactBoolCombStabSets}
Stable formulas relative to $(\mathcal{M}, B)$ are closed under adjunction 
of parameters from $B$ and under boolean combinations. 
\end{fait}
\proof
The cases for adjunction of parameters from $B$ and for negation are clear. 
Hence, for boolean combinations, it suffices to prove the case for disjunction. 
Notice that in the absence of compactness one cannot proceed as in 
\cite[Lemma 0.2.10]{Wagner(Book)97}.

Let $\phi(\overline{x},\overline{y})$ and $\psi(\overline{x},\overline{y})$ 
be two formulas stable relative to $(\mathcal{M}, B)$, 
with stability indices $n_{\phi}$ and $n_{\psi}$ respectively. 
Let $\mu>\max\{n_{\phi},n_{\psi}\}$. 
We claim that the Ramsey number $R(2,2,\mu)$ is a bound for the stability index of 
$[\phi\vee\psi]$ relative to $(\mathcal{M}, B)$.

Assume towards a contradiction there exist tuples 
$\overline{a}_{1}, \ldots , \overline{a}_{i}, \ldots , \overline{a_m}$ and 
$\overline{b}_{1}, \ldots , \overline{b}_{j}, \ldots , \overline{b_m}$ in $B$, 
with $m>R(2,2,\mu)$, such that $[\phi\vee\psi](\overline{a}_{i},\overline{b}_{j})$ 
holds if and only if $i\leq j$. Attach to each pair $\{i,j\}$ from $\{1,\cdots ,m\}$ a 
color, green if $\phi(\overline{a}_{i},\overline{b}_{j})$ or 
$\phi(\overline{a}_{j},\overline{b}_{i})$ holds, and red if 
$\psi(\overline{a}_{i},\overline{b}_{j})$ or 
$\psi(\overline{a}_{j},\overline{b}_{i})$ holds. Notice that, by assumption, 
each pair $\{i,j\}$ has a color (green, red, or both). 
By Fact \ref{finRam}, there exists a subset of 
$\{1,\cdots ,m\}$ of size at least $\mu$ and whose pairs are monochromatic. 
As $\phi(\overline{a}_{i},\overline{b}_{j})$ and 
$\psi(\overline{a}_{i},\overline{b}_{j})$ never hold for $i>j$, we get that on 
the new monochromatic subset 
$\phi(\overline{a}_{i},\overline{b}_{j})$ holds if and only if $i\leq j$ or 
$\psi(\overline{a}_{i},\overline{b}_{j})$ holds if and only if $i\leq j$. But this 
is a contradiction to the fact that $\mu > \max\{n_{\phi},n_{\psi}\}$. 
\qed

\bigskip
We say that a structure $\cal M$ is {\em qf-stable} (``quantifier-free stable") 
if quantifier-free formulas are all stable relative to $(\mathcal{M}, M)$, where $M$ is the 
domain of $\cal M$. This corresponds to the usual notion of stability of $\cal M$ for 
quantifier-free formulas, and by Fact \ref{FactBoolCombStabSets} this is equivalent to the fact that atomic formulas $\phi(\overline{x},\overline{y})$ without 
parameters define stable sets in $\cal M$ in the usual sense. Of course, this is expressed 
in the universal theory of $\cal M$. 

The main ingredient for the proof of Theorem \ref{MainTheo} is the following 
technical lemma. 

\begin{lem}\label{LemInt}
Let $G$ and $H$ be two qf-stable groups and let $w(\overline{x},\overline{y})$ 
be a group word. Assume that, for some separation 
$$\overline{x}=(\overline{x}^{G},\overline{y}^{H}) \mbox{~and~}
\overline{y}=(\overline{y}^{G},\overline{y}^{H}) \leqno(1)$$
of the variables $\overline{x}$ and $\overline{y}$, $w$ has the form 
$$w(\overline{x},\overline{y})=
u_{1}(\overline{x}^{\epsilon_{1}},\overline{y}^{\epsilon_{1}})
\cdots 
u_{k}(\overline{x}^{\epsilon_{k}},\overline{y}^{\epsilon_{k}})
\cdots
u_{\ell}(\overline{x}^{\epsilon_{\ell}},\overline{y}^{\epsilon_{\ell}})
\leqno(2)$$
where $\ell \geq 1$ and, for $1\leq k \leq \ell$, the $\epsilon_{k}$'s 
represent alternatively the symbol $G$ or $H$. 
Then there exists $n$ {\it (depending only on the decomposition of $w$ as 
in $(1)$ and $(2)$)} 
bounding the set of natural numbers $m$ for which there exists natural interpretations 
$$\overline{a}_{i}=(\overline{a_{i}}^{G},\overline{a_{i}}^{H}) 
\mbox{~and~}
\overline{b_{j}}=(\overline{b_{j}}^{G},\overline{b_{j}}^{H})$$ 
in $G\ast H$ of 
$\overline{x}=(\overline{x}^{G},\overline{y}^{H})$ 
and 
$\overline{y}=(\overline{y}^{G},\overline{y}^{H})$ respectively, $1\leq i,j\leq m$, 
such that for each $i$ and $j$ in $\{1,\cdots , m\}$ we have: 
\begin{itemize}
\item
$w(\overline{a_{i}},\overline{b_{j}})=1$ if and only if $i\leq j$, and 
\item
$\overline{a_{i}}^{G},~\overline{b_{j}}^{G} \in G$ and 
$\overline{a_{i}}^{H},~\overline{b_{j}}^{H} \in H$. 
\end{itemize}
\end{lem}
\proof
We proceed by induction on $\ell$. 
For $\ell=1$ everything occurs in a same factor, $G$ or $H$, 
and thus our claim in this case follows from the qf-stability of $G$ and $H$. 

Assume now one has a counterexample $w(\overline{x},\overline{y})$, with 
the corresponding $\ell >1$ minimal. 
By inductive assumption, for each {\em proper} formal sub-expression from the product 
$$w=u_{1}\cdots u_{\ell}$$ 
there is a bound $m$ on the existence of elements with our conditions 
(for this formal sub-expression). 
As in Fact \ref{FactBoolCombStabSets}, there is also such a bound when one considers 
the negation of such formal sub-expressions. In other words, for each proper product 
$\Pi {u_{i}}$ from $u_{1}\cdots u_{\ell}$, and where the factors with same exponent 
$G$ or $H$ are concatenated, we get a bound on sets of elements as in the lemma 
for the formulas $\Pi {u_{i}}=1$ and $\Pi {u_{i}}\neq 1$. 
Let $\mu$ be a (finite) natural number bigger than the maximum of all these bounds. 

We claim that the Ramsey number $R(4\ell, 2, \mu)$ has the desired property. 

Otherwise, one finds 
$\overline{a}_{1}$, ..., $\overline{a_{i}}$,  ..., $\overline{a}_{m}$ and 
$\overline{b}_{1}$, ..., $\overline{b_{j}}$,  ..., $\overline{b}_{m}$, 
with $m>R(4\ell, 2, \mu)$, such that 
$$w(\overline{a_{i}},\overline{b_{j}})=1\mbox{~iff~}i\leq j$$
Now colour the pairs $\{i,j\}$ from $\{1,..., m\}$ with $4\ell$ colours describing, 
for each $k\in \{1,\cdots , \ell\}$, when:  
\begin{itemize}
\item[(A)]
$u_{k}(\overline{a_{i}}^{\epsilon_{k}},\overline{b_{j}}^{\epsilon_{k}})=1$ 
and 
$u_{k}(\overline{a_{j}}^{\epsilon_{k}},\overline{b_{i}}^{\epsilon_{k}})\neq 1$ 
for $i<j$; 
\item[(B)]
$u_{k}(\overline{a_{i}}^{\epsilon_{k}},\overline{b_{j}}^{\epsilon_{k}})\neq 1$ 
and 
$u_{k}(\overline{a_{j}}^{\epsilon_{k}},\overline{b_{i}}^{\epsilon_{k}})=1$ 
for $i<j$; 
\item[(C)]
$u_{k}(\overline{a_{i}}^{\epsilon_{k}},\overline{b_{j}}^{\epsilon_{k}})\neq 1$ 
and 
$u_{k}(\overline{a_{j}}^{\epsilon_{k}},\overline{b_{i}}^{\epsilon_{k}})\neq 1$;  
\item[(D)]
$u_{k}(\overline{a_{i}}^{\epsilon_{k}},\overline{b_{j}}^{\epsilon_{k}})=1$ 
and 
$u_{k}(\overline{a_{j}}^{\epsilon_{k}},\overline{b_{i}}^{\epsilon_{k}})=1$. 
\end{itemize}
Notice that each pair $\{i,j\}$ is well, and uniquely, coloured in this process. By definition of $R(4\ell, 2, \mu)$ and Fact \ref{MainTheo}, 
there exists a monochromatic subset of $\{1,\cdots , m\}$ of size at least $\mu$. 

On each ``coordinate" $u_{k}$, the two first colors (A) and (B) as 
above are excluded. This follows from the fact that $\mu$ is larger then 
the stability indices of the formulas $u_{k}=1$ and $u_{k}\neq 1$ in the 
group $\epsilon_{k}$. This means that on each coordinate $u_{k}$ 
the value of $u_{k}$ is always $1$ or different from $1$ on our monochromatic subset. 

We now claim that at least one coordinate $u_{k}$ is constantly equal to $1$ on 
our monochromatic subset. Otherwise, all interpertations in $G\ast H$ 
of the formal expression 
$$w=\Pi_{k=1}^{\ell}u_{k}$$
would give rise to a normal form in the free product $G\ast H$, 
as all $u_{k}$'s would have nontrivial interpretations alternatively in $G$ or $H$ 
(and as $\ell >1$!). Then one would get a nontrivial value for all 
terms $w(\overline{a_{i}},\overline{b_{j}})$. 
This is a contradiction as approximatively half of them, those for which $i>j$, 
are trivial. 

Now one can discard the coordinates $u_{k}$'s constantly equal to $1$, and 
our induction hypothesis applies to the remaining {\em proper} subword of 
$w=u_{1}\cdots u_{\ell}$. We then get a contradiction, as $\mu$ is larger than 
the stability index (in our situation) of this proper subword. 
\qed

\bigskip
With Lemma \ref{LemInt} one can prove Theorem \ref{MainTheo} as follows. 

\bigskip
\noindent
{\bf Proof of Theorem \ref{MainTheo}.} 
Assume towards a contradiction: For each $m$, there exists 
$\overline{a}_{1}$, ..., $\overline{a}_{i}$, ..., $\overline{a}_{m}$ and 
$\overline{b}_{1}$, ..., $\overline{b}_{j}$, ..., $\overline{b}_{m}$ 
in $B_{r}(G\ast H)$ such that $w(\overline{a_{i}},\overline{b_{j}})=1$ 
if and only if $i\leq j$. 
We will contradict Lemma \ref{LemInt} by making an appropriate 
``change of variables". 

Each element of $B_{r}(G\ast H)$ can be written as a product $r+1$ 
of elements of $G$ and $H$ (alternatively). That is, each element $z$ of 
$B_{r}(G\ast H)$ has the form 
$$z=z_{1}^{G}z_{2}^{H} \cdots z_{r+1}^{G\mbox{~or~}H} \leqno(*)$$ 
with the factors in the product alternatively in $G$ and $H$, as indicated 
by the notation in exponent. 
Now each variable from the tuples $\overline{x}$ and $\overline{y}$
 is replaced by $r+1$ variables as in ($*$) 
 (in particular we multiply by $r+1$ the original number of variables involved in 
 $w(\overline{x},\overline{y})$). 
 Clearly also, each element of $B_{r}(G\ast H)$ has a natural (though not unique) 
 interpretation in $G\ast H$ with this new set of variables, as in equality ($*$). 

Now in the original group word $w(\overline{x},\overline{y})$ we replace formally 
each variable by its expression with the new variables as in equality ($*$), and build 
a new corresponding word $w'$ with these new variables as follows: 
when an inverse of an original variable appears one takes the visual inverse of the 
expression as in ($*$), and for products one just takes concatenation. 
Notice however that one does not proceed to simplifications with the new variables 
as one would do with elements in the group. Then we factor $w'$ as 
$$w'=u_{1}'\cdots u_{\ell}'$$ 
where each subword $u_{k}'$ corresponds to a maximal concatenation of the 
new variables or their inverse with the same exponent $G$ or $H$. That is, 
$u_{1}'$ is the subword corresponding to the first variables with the same exponent, 
$u_{2}'$ is the subword consisting of the next variables with the same exponent, ...etc. 

We are now in position to apply Lemma \ref{LemInt} with this new word $w'$ and 
its decomposition $w'=u_{1}'\cdots u_{\ell}'$ in words $u'_{k}$'s involving alternatively 
variables attached to $G$ or $H$. 
Notice that one gets a contradiction to the statement of Lemma \ref{LemInt} 
with the interpretation of the $\overline{a}_{i}$'s and $\overline{b}_{j}$'s in 
the new variables as in equality ($*$). 

This concludes the proof of Theorem \ref{MainTheo}. 
\qed

\begin{rem}
Unfortunately, the above proof of Theorem \ref{MainTheo} 
highly depends on $r$.  With $w$ fixed, the bound $n$ in Theorem \ref{MainTheo} 
increases as $r$ increases. We note however that $r$ appears only in the change of 
variables in our proof, and in particular not via Lemma \ref{LemInt}. 
Of course, a proof of the qf-stability of the free product of two (or more) 
qf-stable groups requires a technic which does not rely on such a fixed change of variables. 
\end{rem}

We finish with a more general version of Theorem \ref{MainTheo}, where the 
number of factors is not restricted to two. Of course, the notion of length of an element and 
of ball in a free product have natural generalizations for free products of an arbitrary 
number of groups. 

\begin{theo}
Let $\{G_{s}: 1\leq s\leq k\}$ be a family of qf-stable groups and let $\ast G_{s}$ be their free product. Let $w(\overline{x},\overline{y})$ be a group word 
and let $r\geq 1$ be a natural number. Then there exists a natural number $n$ 
(depending on $w$ and $r$) bounding the set of $m$ for which there exists 
$\overline{a}_{1}$, ..., $\overline{a}_{i}$, ..., $\overline{a}_{m}$ and 
$\overline{b}_{1}$, ..., $\overline{b}_{j}$, ..., $\overline{b}_{m}$ 
in $B_{r}(\ast G_{s})$ such that $w(\overline{a_{i}},\overline{b_{j}})=1$ 
if and only if $i\leq j$. 
In particular any quantifier-free formula (possibly with parameters) is stable relative to 
$(\ast G_i, B_{r}(\ast G_i))$. 
\end{theo}

\proof 
One may proceed as in the proof of Theorem \ref{MainTheo}. 
For the change of variables, one now needs to decompose each element of 
$B_{r}(\ast G_i)$ as a product involving variables attached to each $G_{s}$. This may 
give an expression longer than that in equality ($*$) in the proof of 
Theorem \ref{MainTheo}, but is still possible as everything remains finite. 
One then gets a new word $w'$, decomposed 
as $w'=u_{1}'\cdots u_{\ell}'$ where the new variables of each $u'_{k}$ correspond 
alternatively to different factors. One can then prove a similar version of 
Lemma \ref{LemInt}, with possibly a separation of the variables into $k$ blocks 
instead of $2$, similarly by induction on $\ell$. 

Our last statement is just Fact \ref{FactBoolCombStabSets}. 
\qed

\bibliographystyle{alpha}
\bibliography{biblio}

 \textsc{Universit\'{e} Lyon 1, Universit\'{e} de Lyon, CNRS UMR 5208 Institut Camille Jordan, B\^{a}timent du Doyen Jean Braconnier, 43, blvd du 11 novembre 1918, F-69622 Villeurbanne Cedex, France}

\it{azadeh@math.univ-lyon1.fr}

\end{document}